\documentclass[11pt]{article}

\usepackage{amsfonts}
\usepackage{amsmath}
\usepackage{amssymb}

\newtheorem{Theorem}{Theorem}[section]
\newtheorem{Lemma}[Theorem]{Lemma}
\newtheorem{Corollary}[Theorem]{Corollary}
\newtheorem{Proposition}[Theorem]{Proposition}
\newtheorem{Example}[Theorem]{Example}
\newtheorem{Remark}[Theorem]{Remark}
\newtheorem{Conjecture}[Theorem]{Conjecture}
\newtheorem{Question}[Theorem]{Question}

\newenvironment{Proof}{\noindent {Proof:}}{\hfill$\square$\medskip}

\newcommand{\N}{\mathbb{N}}
\newcommand{\Z}{\mathbb{Z}}
\newcommand{\R}{\mathbb{R}}

\newcommand{\OO}{\operatorname{O}}

\newcommand{\C}{\operatorname{C}}

\newcommand{\Irr}{\operatorname{Irr}}
\newcommand{\IBr}{\operatorname{IBr}}
\newcommand{\Hom}{\operatorname{Hom}}

\newcommand{\M}{{\rm M }}

\newcommand{\Alt}{{\rm A }}
\newcommand{\RR}{{\rm R }}
\newcommand{\GL}{{\rm GL}}

\newcommand{\PSL}{{\rm PSL}}
\newcommand{\PGL}{{\rm PGL}}

\newcommand{\I}{{\rm I}}
\newcommand{\Sym}{{\rm S}}

\newcommand{\Nor}{{\rm N}}

\newcommand{\Fix}{{\rm Fix(\it P)}}
\newcommand{\Heart}{{\rm H}}
\newcommand{\Kernel}{{\rm Ker \,}}
\newcommand{\Syl}{{\rm Syl}}

\author{Conchita Mart\'{\i}nez-P\'{e}rez\footnote{Partially supported by BFM2010-19938-C03-03, Gobierno de Arag\'on and European Union's ERDF funds.} \ and Wolfgang Willems\footnote{The author thanks the Department of Mathematics 
          at the University of Zaragoza for its hospitality while this work was completed.}}

\title{On the dimensions of PIM's }

\date{}

\begin{document}
%%%%%%%%%%%%%%%%

\maketitle

\noindent
{\small {\bf Abstract}
We discuss the structure of finite groups for which the projective indecomposable modules have special given dimensions.
In particular, we prove the converse of Fong's dimension formula for $p$-solvable groups. Furthermore,
we characterize  groups for which all irreducible $p$-Brauer characters have $p$-power degrees.


\noindent
{\small {\bf AMS Classification (2010):} 20C20, 20C05 20C33, 20C30}

%%%%%%%%%%%%%%%%%%%%%%%%%%%%%%%%%%%%%%

\section{Introduction}

All  groups in this note are assumed to be finite. The study of character degrees is a relevant subject in character theory. Many of the known results provide characterizations of groups for which a given property of their character degrees holds true. This approach has also been used in modular character theory, although in general the situation is more subtle than in the ordinary case.

Throughout this paper $p$ is always a prime number. 
By $\Irr(G)$, respectively
$\IBr_p(G)$ we denote the irreducible complex, respectively irreducible
$p$-Brauer characters
of the group $G$.
To each $ \varphi \in \IBr_p(G)$ we may associate the projective
indecomposable module  in characteristic $p$ whose head
affords $\varphi$ as Brauer character (called the $\varphi$-PIM). We
use the notation $\Phi_\varphi$ to denote its
complex character.
In particular, $\Phi_\varphi(1)$ is the dimension of the $\varphi$-PIM. 
In case $\varphi $ is the trivial $p$-Brauer character we simply
write $\Phi_1$. Furthermore, we define (after Malle and Weigel, \cite{malleweigel})
$$ c_\varphi = \frac{\Phi_\varphi(1)}{|G|_p} $$
for $ \varphi \in \IBr_p(G)$. In the following let $P$ be a Sylow $p$-subgroup of $G$ and let
$ \varphi^{\Fix}(1)$ denote the dimension of the vector space of $P$-fixed points of
the module afforded by $\varphi$. Obviously,
\begin{equation} \label{fix} 1 \leq  \varphi^{\Fix}(1) \leq (1, \Phi_\varphi|_P) = (1_P^G,\Phi_\varphi) = c_\varphi. \end{equation}

 For $c_\varphi$, the following three situations are of particular interest. 
$$ \begin{array}{ll}
 {\rm (i)} \ \ \qquad c_\varphi = \varphi(1) & \qquad \mbox{for all} \ \varphi \in \IBr_p(G), \\[1ex]
 {\rm (ii)} \  \qquad c_\varphi = \varphi(1)_{p'} & \qquad \mbox{for all} \ \varphi \in \IBr_p(G), \\[1ex]
 {\rm (iii)} \qquad c_\varphi = \varphi^\Fix(1) & \qquad \mbox{for all} \ \varphi \in \IBr_p(G), \\[1ex]
 {\rm (iv)} \ \qquad c_\varphi = 1  & \qquad \mbox{for all} \ \varphi \in \IBr_p(G).
\end{array}
$$
The first case  has already been considered  by Brockhaus in \cite{Brockhaus}. He proved that 
$c_\varphi = \varphi(1) $ for all $  \varphi \in \IBr_p(G)$ if and only if the Sylow $p$-subgroup of
$G$ is normal. In \cite{Willems} we conjectured 
$$   |G|_{p'} \leq \sum_{\varphi \in \IBr_p(G)} \varphi(1)^2 $$
with equality if and only if the Sylow $p$-subroup is normal. Unfortunately this
conjecture is still open and we will discuss it slightly more in the last section of this paper.

\bigskip

Clearly, the second condition (ii) holds true for $p$-solvable groups, by Fong's dimension formula (see
\cite{Navarro}, Corollary 10.14). To be brief let $\IBr_p(\Phi_1)$ denote the set of irreducible $p$-Brauer characters which occur as constituents in $\Phi_1$. 
We prove in section \ref{p'}

\bigskip

\noindent
{\bf Theorem A} 
If the finite group $G$ satisfies
\begin{equation} \label{p-solv}  c_\varphi = \varphi(1)_{p'} \quad
\mbox{for all} \  \varphi \in \IBr_p(\Phi_1)  \end{equation}  then $G$ is $p$-solvable.
In particular, condition (ii) characterizes  $p$-solvable groups.
\bigskip

At a first glance  condition (iii) does not seem to be very natural. However, as one of the main results we prove in section \ref{Fix(P)}

\bigskip
\noindent
{\bf Theorem B}  If the finite group $G$ satisfies
\begin{equation} \label{p-solv2}  c_\varphi = \varphi^\Fix(1) \ \mbox{for all} \  \varphi \in \IBr_p(\Phi_1)
\end{equation}
then $G$ is $p$-solvable. Moreover, $G$ satisfies (\ref{p-solv2}) for every $\varphi\in\IBr_p(G)$ if and only if the permutation module $1_P^G$  is
completely reducible  in characteristic $p$. 

\bigskip
Finally, in case (iv), the condition $c_\varphi = 1$ implies $\varphi^\Fix(1) = c_\varphi = 1$, by 
(\ref{fix}). Thus, if this condition holds for all $\varphi \in \IBr_p(G)$  then $G$ is $p$-solvable according to Theorem B,
and Fong's Dimension Formula yields
$$ 1 = c_\varphi = \varphi(1)_{p'}.$$
This means that all irreducible $p$-Brauer characters have   $p$-power degrees.
If $p$ is odd then $G$ has to be solvable, by Corollary 2 of \cite{Tiep-Willems}.
Since a $2$-solvable group is solvable by the Feit-Thompson Theorem, we see that
$G$ is a solvable group. Such groups with $c_\varphi = 1$ for  all $\varphi \in \IBr_p(G)$
have been considered already in section 13 of \cite{ManzWolf}. In particular, in Corollary 13.10 it is proved that
$G=\OO_{p,p',p,p',p}(G)$.  We give the full classification in section \ref{trivial},
namely

\bigskip
\noindent
{\bf Theorem C} The following conditions are equivalent.
\begin{itemize}
\item[a)] $ c_\varphi = 1$ for all $\varphi \in \IBr_p(G)$.
\item[b)] $ G=\OO_{p,p',p,p',p}(G)$ where $L=\OO_{p,p',p,p'}(G)/\OO_p(G)$ satisfies the following two conditions.  
                \begin{itemize}
                        \item[$\alpha$)] All $p'$-sections of $L$ are abelian.

                         \item[$\beta)$]  $L=C_L(x)\OO_{p',p}(L)$ for all $x\in \OO_{p'}(L)$.
                \end{itemize}
\end{itemize}
 
In the literature there are several results which prove that properties of the classical irreducible constituents
of $\Phi_1$ heavily restrict the structure of the underlying group; see for instance \cite{Pahlings}, \cite{MW}.
Theorem A and Theorem B show that properties on all irreducible Brauer characters occuring in $\Phi_1$ also have a strong impact on the 
group structure.

\section{The case $c_\varphi = \varphi(1)_{p'}$} \label{p'}

We start this section with a useful lemma.

\begin{Lemma}\label{normal} Let $N$ be a normal subgroup of $G$. 
\begin{itemize}
\item[\rm a)]  If $\varphi\in\IBr_p(G)$ is obtained by inflation from $\bar\varphi\in\IBr_p(G/N)$ then $c_\varphi=c_{\bar\varphi}\,c_{1_N}$.
\item[b)] Let $\theta\in\IBr_pN$ and $\varphi\in\IBr_pG$ such that $\theta$ is a constituent of $\varphi_N$. Then
$$   c_\varphi = \left\{ \begin{array}{ll}
                   {\varphi(1)\over\theta(1)}c_\theta & \mbox{if} \ p\nmid|G:N| \\
                      c_\theta & \mbox{if} \ p=|G:N|. 
         \end{array} \right.             
$$
\end{itemize}
\end{Lemma}
\begin{Proof} 
a)  By (\cite{Huppert-Black}, Chap. VII, 14.2), we have $\Phi_\varphi(1) = \Phi_{\overline{\varphi}}(1)\, \Phi_{1_N}(1)$
where $\Phi_{1_N}$ denotes the character of the $1$-PIM of $N$. \\
b) Observe that the inertial groups $I=I_G(\theta)$ and $I_G(\Phi_\theta)$ coincide.  If $p\nmid|G:N|$ 
then by
  (\cite{Navarro}, Corollary 8.7 and Corollary 8.8), we immediately get $\Phi_\varphi(1)={\varphi(1)\over\theta(1)}\,\Phi_\theta(1).$
In case $|G:N|=p$  it suffices to note that Green's Indecomposability Theorem (\cite{Huppert-Black}, Chap. VII,
16.2) yields $\Phi_\theta^G=\Phi_\varphi$.
\end{Proof}

\begin{Lemma} \label{factorgroup}
Let $N$ be a normal subgroup of $G$. If $G$ satisfies
$$  c_\varphi =\varphi(1)_{p'} $$
for all $ \varphi \in \IBr_p(\Phi_1)$  then
  the factor group $G/N$ satisfies the analogous equations. If moreover $G/N$ is $p$-solvable, then so does $N$.
\end{Lemma}
\begin{Proof} Note that the hypothesis implies $c_1=c_{1_N}=1$.
Let $\Phi_{\overline{1}}$ denote the
character of the $1$-PIM of $\overline{G} = G/N$ and let
$\overline{\varphi} \in \IBr_p(\overline{G})$ be in
$\IBr_p(\Phi_{\overline{1}})$.
If $\varphi$ is the inflation of $\overline{\varphi}$ to $G$ then
obviously $\varphi \in \IBr(\Phi_1)$. Thus the first assertion follows by Lemma \ref{normal} a).

For the last part of the Lemma, let $\theta$ be a constituent of $\Phi_{1_N}$. Since $G/N$ is $p$-solvable by assumption, we may assume that either $p\nmid|G:N|$ or  $|G:N|=p$. We consider  first the case $ p \nmid |G:N|$. As  $\Phi_1|_N=\Phi_{1_N}$  there exists some $\varphi\in \IBr_p(\Phi_1)$ where $\theta$ is a constituent of $\varphi_N$ (this also follows  from Corollary 8.7 and Corollary 8.8 in \cite{Navarro}, or from \cite{Huppert-Black}, Chap. VII, Lemma 14.2).
According to the assumption of the Lemma  we have $c_\varphi=\varphi(1)_{p'}$. 
 Moreover, Dade's Theorem (\cite{Navarro} Theorem 8.30) says that $p\nmid {\varphi(1)\over\theta(1)}$. 
Thus Lemma \ref{normal} b) yields $c_\theta=\theta(1)_{p'}$.

In case $|G:N|=p$
 Green's Indecomposability Theorem (\cite{Huppert-Black}, Chap. VII,
16.2) yields $\Phi_{1_N}^G=\Phi_1$.  Thus we  may choose again a $\varphi\in \IBr_p(\Phi_1)$ where $\theta$ is a constituent of $\varphi_N$. In this case, either $\varphi(1)=\theta(1)$ or $\varphi(1)=p\theta(1)$, hence $\varphi(1)_{p'} = \theta(1)_{p'}$, and the assertion follows by Lemma \ref{normal} b). 
\end{Proof}

Suppose that the next result has been proved already.

\begin{Proposition} \label{simple} Let $G$ be a simple non-abelian group. If
$$  c_\varphi = \varphi(1)_{p'} $$
for all $ \varphi \in \IBr_p(\Phi_1)$  then $G$ is a $p'$-group.
\end{Proposition}

\bigskip

\noindent
{\bf Proof of Theorem A} Let $G$ be a non $p$-solvable group
of minimal order which satisfies the assumptions of Theorem A.
By Proposition \ref{simple}, the group $G$ is not simple. Now we take any  nontrivial normal subgroup $N$ of $G$. By Lemma \ref{factorgroup}, the factor group
$G/N$ is $p$-solvable. Under this hypothesis Lemma \ref{factorgroup} implies that $N$ is  $p$-solvable as well which completes the proof.

\bigskip
Thus we are left with the proof of Proposition \ref{simple}.
 Let $G$ be a simple non-abelian group satisfying
\begin{equation} \label{degree}  c_\varphi = \varphi(1)_{p'} \text{ for all $ \varphi \in \IBr_p(\Phi_1)$.}\end{equation}
If  $ p \mid |G|$ then, by Theorem A of \cite{malleweigel}, the group
$G$ belongs to one in the following list which we call  the Malle-Weigel list in the sequel.  \\
\begin{itemize}
  \item[(a)] $\Alt_p$, $ p \geq 5$
  \item[(b)] $\PSL(2,p)$, $ p \geq 5$
  \item[(c)] $\PSL(n,q)$ where $\frac{q^n-1}{q-1} = p^f$ and $n \geq 3$
  \item[(d)] $\PSL(2,q)$ where $q$ is a Mersenne prime and $p=2$
  \item[(e)] $\PSL(2,q)$ where $p=q+1$ is a Fermat prime
  \item[(f)] $\PSL(2,8)$ and $p=3$
  \item[(g)] $\M_{11}$, $p=11$
  \item[(h)] $\M_{23}$, $p=23$
\end{itemize}

In each case  we will find a $ \varphi \in \IBr_p(\Phi_1)$ which
does not satisfy
the equation (\ref{degree}). This can  be seen directly for
$\M_{11}, \M_{23}$ and $\PSL(2,8)$ using the GAP library \cite{GAP}.
It also follows immediately in case (e) using (\cite{Burkhardt}, case
II) and in case (d) using (\cite{Burkhardt}, case VIII (b)).
Furthermore, $G=\PSL(2,p)$, $p \geq 5$ does not satisfies (\ref{degree}) either, since by (\cite{Alperin}, Section 7),
there is a $\varphi \in \IBr_p(\Phi_1)$ with $\varphi(1) = p-2 = \varphi(1)_{p'} > 2$ and
$$ c_\varphi = \frac{\Phi_\varphi(1)}{p} = 2 \not= \varphi(1)_{p'}.$$

Thus we only have to deal with the cases in (a)  and (c). But for further reference, we shall consider the whole list in the next result.

\begin{Lemma} \label{pim} For each of the groups $G$ in the Malle-Weigel list we have the following.
\begin{itemize}
   \item[\rm a)] The Sylow $p$-subgroups of $G$ are cyclic except in case {\rm (d)}.
   %\item[b)] $\Phi_1(1) = |G|_p.$
   \item[\rm b)] $\Phi_1$ contains exactly one non-trivial Brauer
character $\beta$ with $\beta(1) = |G|_p -2 = \beta(1)_{p'}$
             except possibly in case {\rm (d)}.
\end{itemize}
\end{Lemma}
\begin{Proof} a) This is clear for $\Alt_p$ and $\PSL(2,p)$. The
condition $\frac{q^n-1}{q-1} = p^f$ for $n \geq 3$ in case (c) implies, by
Zsigmondy's Theorem (see \cite{Huppert-Black}, Chap. IX, Theorem 8.3), that
$|\PSL(n,q)|_p = p^f$ since $n \geq3$ and $2^6-1 \not= p^f$.
Moreover, in this case a Sylow $p$-subgroup is generated by a
Singer cycle. Finally, for all groups in (e) - (h) the assertion is clear.  \\
%b) This is obvious, since $c_1 = 1$. \\
b) For $\Alt_p$ it follows by a result of Wielandt (see
\cite{Huppert-Black}, Chap. XII, Theorem 10.7).
In case (b) and (c) we can read it off from the Brauer trees given in
(\cite{Alperin}, p. 123) and (\cite{FS}, Theorem C).
Finally, case (e) and (f) follow by \cite{Burkhardt}, case (g) and (h) by
the GAP library \cite{GAP}.
\end{Proof}

\begin{Lemma} \label{general3} Assume that $\Phi_1(1)=|G|_p$. If there is  a $\beta \in \IBr_p(G)$
with $\beta\not\in\IBr_p(\Phi_1)$ such that
$\Phi_\beta\mid \Phi_1\Phi_1$, then
for some $\varphi \in \IBr_p(\Phi_1)$
we have $c_\varphi < \varphi(1)$. \\
In particular, {\rm (\ref{degree})} does not hold true if in addition  $p
\nmid \varphi(1)$ for every $\varphi \in \IBr_p(\Phi_1)$.
\end{Lemma}
\begin{Proof} As for any $\varphi\in\IBr_p(G)$ the character $\Phi_\varphi$ is a constituent of the projective character $\varphi\Phi_1$, the first assumption implies $c_\varphi \leq \varphi(1)$.
Suppose that $c_\varphi = \varphi(1)$ for all $\varphi \in
\IBr_p(\Phi_1)$. Thus $\Phi_\varphi = \varphi \Phi_1$.
So, if we put $ \Phi_1 = \sum_{\varphi \in \IBr_p(\Phi_1)}
a_\varphi \varphi $, then
$$ \Phi_1  \Phi_1 = \sum_{\varphi \in \IBr_p(\Phi_1)} a_\varphi \
\varphi  \Phi_1 =
         \sum_{\varphi \in \IBr_p(\Phi_1)} a_\varphi  \ \Phi_\varphi,$$
which contradicts the hypothesis.
\end{Proof}

Note that the assertion of the previous Lemma holds true if
$\Phi_1\Phi_1$ contains an irreducible character
of $p$-defect zero.

\begin{Lemma} \label{alternating} Let $p \geq 5$ be a prime.
For $\Alt_p$ the character $\Phi_1 \Phi_1$ contains an irreducible constituent
which belongs to a $p$-block of defect zero.
In particular, $\Alt_p$ does not satisfy {\rm (\ref{degree})}.
\end{Lemma}
\begin{Proof} Let $G= \Sym_p$ be the symmetric group on $p$ letters and let $H=\Sym_{p-1}$. Recall that the classical irreducible characters of $G$ are
labeled by the partitions of $p$.
By the hook length formula (see \cite{FH}, 4.12), the irreducible character
$\phi_{(p-2,2)}$ (corresponding to the partition $(p-2,2)$) belongs to a p-block
of defect zero. Furthermore $1_H^G = \Phi_1$. Using Pieri's
formula (see \cite{FH}, part I, section 4) and Frobenius reciprocity we obtain
$$ \begin{array}{rcl}
1_H^G \, 1_H^G  & = & (1_H^G|_H)^G = ((\phi_{(p)} +
\phi_{(p-1,1)})|_H)^G = (2 \phi_{(p-1)} + \phi_{(p-2,1)})^G \\[2ex]
  & = & 2 \phi_{(p)} + 3 \phi_{(p-1,1)} + \phi_{(p-2,2)}+\phi_{(p-2,1,1)}.
\end{array} $$
Thus $\phi_{(p-2,2)}$ is a constituent of $\Phi_1 \Phi_1$ and we may apply
Lemma \ref{general3}. % we get $ c_\beta < \beta(1) = \beta(1)_{p'}$.
Since $\Phi_1|_{\Alt_p}$ is equal to $\Phi_{1_{\Alt_p}}$ and $\phi_{(p-2,2)}$
restricts irreducibly to $\Alt_p$
the proof is complete according to Lemma \ref{pim} b).
\end{Proof}

\begin{Lemma} \label{reduction} Let $ N \leq M \leq G$ where $N,M
\lhd G$ and $H=M/N$. Suppose that $ p \nmid |N|\,|G/M|$. If $\Phi_{1_G}
\Phi_{1_G}$
contains an irreducible character of $p$-defect zero then $\Phi_{1_H}
\Phi_{1_H}$ does as well.
\end{Lemma}
\begin{Proof} According to  (\cite{Huppert-Black}, Chap. VII, Lemma 14.2) we
have $ \Phi_{1_M} = \Phi_{1_G}|_M$. The same reference implies  $\Phi_{1_H} = \Phi_{1_M}$.
Finally observe that the irreducible constituents of the restriction of an
irreducible character of $p$-defect zero to a normal subgroup are also 
 of $p$-defect zero  which completes the proof.
\end{Proof}

\begin{Lemma}\label{generallinear}Let $G=\operatorname{PSL}(n,q)$
with $p^f=(q^n-1)/(q-1)$ and $n \geq 3$. Then $\Phi_1\Phi_1$ has an irreducible
constituent which belongs to a $p$-block of defect zero. In
particular, $G$  does not satisfy {\rm (\ref{degree})}.
\end{Lemma}
\begin{Proof} By Lemma \ref{reduction}, we may assume that $G =
\GL(n,q)$. As observed in the proof of Lemma \ref{pim}, the condition $(q^n-1)/(q-1)=p^f$ implies 
$p^f=|G|_p$ (recall that $n\neq 2$). This means that the parabolic subgroup
$$H=\Bigg\{
\begin{pmatrix}
& a&*\\
&0&\operatorname{GL}(n-1,q)\\
\end{pmatrix} \mid 0 \not= a \in \operatorname{GF}(q)
\Bigg\}$$
is a $p$-complement. Thus $\Phi_1=1_H^G=R^G_L(1)$ where $L$ is the
subparabolic subgroup of $G$ isomorphic to $\operatorname{GL}(1,q)\times\operatorname{GL}(n-1,q)$ and
$\RR^G_L$ denotes the Harish-Chandra induction; i.e.,
$\RR^G_L(\phi)=(\text{infl}\,\phi)^G$ where $\text{infl}\,\phi$ is the
inflation of $\phi \in \Irr(L)$ to $H$. Observe also that the
associated Weyl subgroup $W_H$ is $\Sym_{n-1}$. For any $x\in G\setminus H$
the Mackey formula yields
$$\Phi_1\Phi_1=1_H^G1_H^G=((1_H^G)|_H)^G=1_H^G+1_{H\cap H^x}^G$$
since the action of $G$ on $G/H$ is
2-transitive. If
$$ x=\begin{pmatrix}
&0&1&0\\
&1&0&0\\
&0&0&I_{n-2}\\
\end{pmatrix}$$
then
$$H_1 =H\cap H^x=\Bigg\{
\begin{pmatrix}
& a&0&*\\
&0& b &*\\
&0&0&\GL(n-2,q)\\
\end{pmatrix} \mid 0 \not= a, b \in \operatorname{GF}(q)
\Bigg\}. $$
Assume first that $n>3$. The group $H_1$ is contained in the parabolic subgroup associated to the
partition $\lambda=(n-2,2)$, i.e.,

$$H_1\leq P_{\lambda}=\Big\{
\begin{pmatrix}
&\GL(2,q)&*\\
&0&\GL(n-2,q)\\
\end{pmatrix}
\Big\}.
$$
Thus $ 1_{H_1}^G = 1_{P_{\lambda}}^G + \rho$ for some character $\rho$. Let
$\chi_{\lambda}$ be the unipotent character of $G$ associated to the partition
$\lambda$. As $1_{P_{\lambda}}^G=\RR_{L_{\lambda}}^G(1)$ is also a
character obtained by Harish-Chandra induction from the corresponding
subparabolic subgroup $L_{\lambda}$ we deduce, by (\cite{FS}, Lemma A),
that
$$(\chi_{\lambda},1_{P_{\lambda}}^G)=(\chi_{\lambda},\RR_{L_{\lambda}}^G(1))=(\phi_{\lambda},1_{\Sym_{\lambda}}^{\Sym_n})\neq
0$$
where  $\phi_{\lambda} \in \Irr(\Sym_n)$ is the character
corresponding to $\lambda$ and $\Sym_{\lambda}\cong \Sym_2\times S_{n-2}$ is
the associated Young subgroup. Moreover, by the hook length formula
for degrees of unipotent characters (see for example \cite{Green}),
we have 
$$\chi_{\lambda}(1)=q^{d}{\frac{(q^n-1)(q^{n-3}-1)}{(q-1)(q^2-1)}}$$
for a suitable $d$. The condition $(q^n-1)/(q-1) = p^f$ implies that $n$ is odd.
It follows  $q^2-1\mid q^{n-3}-1$ for $ n \geq 5$ and therefore
$p^f\mid\chi_{\lambda}(1)$. Thus the character $\chi_{\lambda}$
belongs to a  $p$-block of defect zero.

We are left with the case $n=3$. We show again that $1_{H_1}^G$ contains an irreducible character of $p$-defect zero.
 The claim may be
proved using the character table of $\GL(3,q)$, which was first
computed by Steinberg in (\cite{Steinberg}, section 3). To simplify
calculations one can proceed as follows: First check that
$$(1_{H_1}^G,\chi_{(2,1)})_G=(1_{H_1},\chi_{(2,1)})_{H_1}=3$$ and
$$(1_{H_1}^G,\chi_{(1,1,1)})_G=(1_{H_1},\chi_{(1,1,1)})_{H_1}=2.$$
As obviously $(1_{H_1}^G,1_G)_G=1$, and
$1_{H_1}^G$ has degree $q(q+1)(q^2+q+1)$ we are left with a character of degree
$$q(q+1)(q^2+q+1)-1-3q(q+1)-2q^3=(q^2-q-1)(q^2+q+1).$$
Moreover, the condition that $q^2+q+1$ is a prime power implies that
$q\not\equiv 1 \bmod 3$ since otherwise  $q^2+q+1 \equiv 0 \bmod 3$, hence
$q^2+q+1=3^f$ and $3\mid q-1$ which is impossible by Zsigmondy's Theorem.
Now, consider any character of $G$ of the form $\alpha\chi_\mu$ where
$1\neq\alpha$ is linear and $\chi_\mu$ is unipotent. Since $Z ={\rm Z}(G)\leq
H_1$ we have $Z\leq\text{Ker}\, \chi_\mu$. %(this can of course also
%be read off the character table).
Furthermore note that $\alpha$ is a character of the form
$\alpha(x)=(\text{det}\, x)^j$ for some $j$. So, if $3j \not\equiv 0 \bmod q-1$
then $\alpha|_Z\neq 1_Z$. Therefore, in this case, we have
$$(1_{H_1}^G,\alpha\chi_{\mu})_G=(1_{H_1},\alpha\chi_{\mu})_{H_1}\leq
(1_Z,\alpha\chi_\mu)_Z=\chi_\mu(1)(1_Z,\alpha)_Z=0.$$
 %(however, one can
%check with a bit longer calculations that also in the general case
%$(1_{H_1}^G,\alpha\chi_{\mu})_G=0$).
This implies that the irreducible constituents 
left of $1_{H_1}^G$ must have
degree  either 
 $(q-1)^2(q+1)$ or a multiple of $q^2+q+1$  (see
(\cite{Steinberg} section 3).
% (these are the cuspidal characters).
Note that $(q-1)^2(q+1)\equiv3$ mod $q^2+q+1$. So if all the rest irreducible components of $1_{H_1}^G$ were of degree $(q-1)^2(q+1)$, as as $q^2+q+1$ is a prime power we would deduce that $q^2-q-1$ is a multiple of $(q-1)^2(q+1)$ which is impossible. 
Thus $1_{H_1}^G$ and
therefore $\Phi_1\Phi_1$ must contain some irreducible character of degree a
multiple of $q^2+q+1$.
  All such characters belong to $p$-blocks of defect zero. Thus we have
shown that in all cases $\Phi_1\Phi_1$ contains an irreducible
  character of $p$-defect zero. 
  Again, the result follows by Lemma \ref{pim} b). \end{Proof}

\section{The case $c_\varphi = \varphi^\Fix(1)$} \label{Fix(P)}

Let $F$ denote the underlying field of characteristic $p$ (always large enough) and
let $P$ be a Sylow $p$-subgroup of $G$. By  abuse of notation we denote by $1_P$ the trivial $FP$-module (and also its ordinary character). For $\varphi \in \IBr_p(G)$ and $V$ the $FG$-module affording $\varphi$,
we put
$$\varphi^\Fix(1)=\dim_F \, \{ v\in V \mid \, vx = v \ \mbox{for all} \ x \in P\}=\dim_F\text{Hom}_{FP}(1_P,V).$$
Observe that this definition does not depend on the chosen $p$-Sylow subgroup. 
The second part of Theorem B is easy to see.

\begin{Proposition} \label{lower} The following conditions are equivalent.
\begin{itemize}
 \item[\rm a)] $c_\varphi = \varphi^\Fix(1)$ for all $\varphi \in \IBr_p(G)$.
 \item[\rm b)] $1_P^G$ is completely reducible.
\end{itemize}  
\end{Proposition}
\begin{Proof} 
Observe that the multiplicity of $\varphi$ in $1_P^G$ is equal to
$$(1_P^G,\Phi_\varphi) = (1_P, {\Phi_\varphi}|_P) =  c_\varphi.$$
On the other hand, the degree $\varphi^\Fix(1)$ equals to the  multiplicity of
$\varphi$ in the head of $1_P^G$, by
Nakayama's Theorem. Thus $c_\varphi = \varphi^\Fix(1)$ for all $\varphi \in \IBr_p(G)$ if and only if
$1_P^G$ is com\-pletely
reducible.
\end{Proof}

We may naturally  ask: Which are exactly the groups that satisfy $c_\varphi = \varphi^\Fix(1)$ for all $\varphi \in \IBr_p(G)$?
%First, we observe that $\varphi^\Fix(1)$ does not change if we factor out normal subgroups in the kernel of $\varphi$.

%\begin{Lemma}\label{kernel} Let $G$ be a group and let $N$ be a normal subgroup of $G$. For every $G/N$-module $V$ we have
%$$\Hom_G(V,1_P^G)=\Hom_{G/N}(V,1_{PN/N}^{G/N}).$$
%\end{Lemma}
%\begin{Proof} As $N$ acts trivially on $V$, Nakayama's rule yields
%$$\Hom_G(V,1_P^G)=\Hom_P(V,1_P)=\Hom_{PN/N}(V,1_{PN/N})=\Hom_{G/N}(V,1_{PN/N}^{G/N}).$$
%\end{Proof}

\begin{Lemma}\label{red1} Let $G$ be a group and let $N$ be a normal subgroup of $G$.
\begin{itemize}
\item[\rm a)]
           If $1_P^G$ is completely reducible then $1_{PN/N}^{G/N}$ and $1_{P\cap N}^N$ are completely reducible as well.
\item[\rm b)] The converse of a) holds true if
     $N$ is a $p$-group or
     $p \nmid |G/N|$.
\end{itemize}
\end{Lemma}
\begin{Proof} a) Since $1_{PN/N}^{G/N}$ is a factor module of $1_P^G$  and  $1_{P\cap N}^N$ is a direct summand of the restriction $1_P^G|_N$ the statement follows. \\
b) If $N$ is a normal $p$-subgroup we have $1_P^G \cong 1_{P/N}^{G/N}$.
If $p$ does not divide $|G/N|$ the assertion  is a consequence of the fact that an $FG$-module $U$ is completely reducible if and only if the restriction $U_N$ is completely reducible (see \cite{Huppert-Black}, Chap. VII, Theorem 7.21).
\end{Proof}

\begin{Lemma}\label{normal2} Let $N$ be a normal subgroup of $G$.
\begin{itemize}
\item[\rm a)]  If $\varphi\in\IBr_p(G)$ is obtained by inflation from $\bar\varphi \in \IBr_p(G/N)$ then $$\varphi^{\Fix}(1)={\bar\varphi}^{{\rm Fix(\it PN/N)}}(1).$$

\item[\rm b)] Let $\varphi\in\IBr_p(G)$. If $\theta\in\IBr_p(N)$ is a constituent of $\varphi_N$  then
$$ \begin{array}{ll}
     \varphi^{\Fix}(1) =  {\varphi(1)\over\theta(1)}\theta^{{\rm Fix(\it P)}}(1) & \quad \mbox{if} \ p\nmid|G:N|  \end{array} $$
 and
$$ \begin{array}{ll}       \varphi^{\Fix}(1)\leq\theta^{{\rm Fix(\it P\cap N)}}(1) & \quad \mbox{if} \ p=|G:N|. 
   \end{array}
$$
\end{itemize}
\end{Lemma}
\begin{Proof} a) This is obvious. \\
b) Let $V$ and $W$ be modules over $F$ affording $\varphi$ and $\theta$ respectively. We first consider the case $p\nmid|G:N|$, hence $P \leq N$. Note that any $G$-conjugate character $\theta^x$ of $\theta$ satisfies $(\theta^x)^{{\rm Fix(\it P)}}(1)=\theta^{{\rm Fix(\it P)}}(1)$. Therefore, by Nakayama's Lemma,
$$\varphi^\Fix(1)=\dim_F\text{Hom}_{FG}(1_P^G,V)=\dim_F\text{Hom}_{FN}(1_P^N,V_N)={\varphi(1)\over\theta(1)}\theta^{{\rm Fix(\it P)}}(1).$$
In  case $p=|G:N|$ we have $G=PN$. Thus, by Nakayama's Lemma and Mackey's Formula, we have 
$$\varphi^\Fix(1)\leq\dim_F\text{Hom}_{FG}(1_P^G,W^G)=\dim_F\text{Hom}_{FN}(1_{P\cap N}^N,W)
=\theta^{{\rm Fix(\it P\cap N)}}(1).$$

\end{Proof}

The next result proves the first part of Theorem B. 

\begin{Theorem} \label{p-solvable} If $\varphi^\Fix(1)=c_\varphi$ for every $\varphi\in\IBr_p(\Phi_1)$ then $G$ is $p$-solvable.
\end{Theorem}
\begin{Proof} 
Let $G$ be a minimal counterexample and note hat $c_1=1^\Fix(1)=1$. We assume first that $G$
 is not simple. Let $1 \not= N \lhd G$. Then Lemma \ref{normal} a) together with Lemma \ref{normal2} a) imply that $G/N$ satisfies the same hypothesis as $G$. Thus we may assume by induction that $G/N$ is $p$-solvable. In particular, we may assume that either $p\nmid|G:N|$ or $p=|G:N|$. If we prove that $N$  satisfies the hypothesis as well then we are done by induction. Let $\theta\in\IBr_p(\Phi_{1_N})$. By the same argument as in Lemma \ref{factorgroup} we see 
in both cases that there is some $\varphi\in\IBr_p(\Phi_1)$ such that $\theta$ is a constituent of $\varphi_N$. The hypothesis of the Theorem says $\varphi^\Fix(1)=c_\varphi$. If $p\nmid|G:N|$ then Lemma \ref{normal} b) together with Lemma \ref{normal2} b) imply the same for $\theta$. In case $|G:N|=p$, Lemma \ref{normal2} b) yields
$$c_\varphi=\varphi^{\Fix}(1)\leq\theta^{{\rm Fix(\it P\cap N)}}(1)\leq c_\theta.$$
According to Lemma \ref{normal} b) we have $c_\varphi=c_\theta$ and we are done.

So we may assume that $G$ is a simple non-abelian group of order divisible by $p$.
Then $G$ is one of the groups in Malle-Weigel's list of Section 3.

First we consider any group in this list but not being in (d). By Lemma \ref{pim}, the Sylow $p$-subgroups 
of $G$ are cyclic.
Let  $Q$ be an $FG$-module affording $\Phi_1$. Clearly, $Q|_P=FP$  since $\Phi_1(1)=|P|$. Furthermore, $Q|_P$ is uniserial
since $P$ is cyclic. It follows that the module affording the only non-trivial irreducible character $\beta\in\Phi_1$ (see Lemma \ref{pim}) is uniserial as well 
when restricted to $P$
which proves that
$c_\beta = \beta^\Fix(1)=1$. Since $\beta(1)=|P|-2$ we obtain
$$|P|=\Phi_\beta(1)=2\beta(1)+m=2|P|-4+m$$
 where $m>1$ is the degree of some Brauer character. Thus $|P| =4-m\leq 2$, which is impossible.

 It remains to deal with the case $G=\PSL(2,q)$ where $q$ is a Mersenne prime and $p=2$.
In particular, $q+1 = 2^t$ for some $t \geq 2$.  If $Q$ denotes again the $FG$-module affording $\Phi_1$
then the heart $\Heart(Q)$ of $Q$ satisfies $\Heart(Q) = V \oplus W$ with irreducible modules $V,W$, both of dimension
$\frac{2^t-2}{2} = 2^{t-1}-1$ (see \cite{Erdmann}, Theorem 4).  Since $Q|_P=FP$ and  $\Heart(FP)$ is a direct sum of two uniserial modules
(note that $P$ is a dihedral group) we see that both modules $V$ and $W$ have a one-dimensional fixed point subspace for $P$.
Thus, by the hypothesis of the Theorem, $\Phi_\varphi(1)=\Phi_\beta(1)=|G|_2=2^t$ where $\varphi$ and $\beta$ are the Brauer characters of $V$ and $W$, respectively. Furthermore, since $\Phi_\varphi(1)\geq 2\varphi(1)=2^t-2$ we see that $\Phi_\varphi$ contains $1_G$ exactly with multiplicity 2.
But, by symmetry of the
Cartan matrix,  this multiplicity equals the multiplicity of $\varphi$ in $\Phi_1$  which is $1$. Thus we have a contradiction which completes the proof.

\end{Proof}

Note that according to Proposition \ref{lower}  the previous result  implies that  $G$ is $p$-solvable if $1_P^G$ is completely reducible.

\begin{Proposition} \label{solvable} Suppose that  $c_\varphi = \varphi^\Fix(1)$ for all $\varphi \in \IBr_p(G)$. If  $H=\OO^{p'}(G)$ then
\begin{itemize}
 \item[\rm a)]  $ p \mid \varphi(1)$ for all $1 \not= \varphi  \in \IBr_p(H)$.
 \item[\rm b)] $\Nor_H(P) =P$.
 \item[\rm c)] $H$ is solvable.
\end{itemize}
\end{Proposition}
\begin{Proof} The hypothesis of the Proposition on $G$ passes on to $H$, by 
Proposition \ref{lower} and Lemma \ref{red1} b).
By Theorem \ref{p-solvable}, the group $H$ is $p$-solvable, and by Lemma \ref{red1}, the module $1_P^H$ is completely reducible. \\
a) Suppose there is some $1 \not= \varphi  \in \IBr_p(H)$ with $ p \nmid \varphi(1)$. Since $H$ is $p$-solvable we obtain
$$ \varphi(1) = \varphi(1)_{p'} = c_\varphi = \varphi^\Fix(1).$$
Thus $P$ is contained in the kernel of $\varphi$ which is therefore a proper normal subgroup of $H$ of index prime to $p$.
This contradicts the fact that $H = \OO^{p'}(G)$. \\
b) In a $p$-solvable group $X$ the number of irreducible $p$-Brauer characters of $p'$-degree equals the number
of irreducible classical characters of $\Nor_X(P)/P$ where $P$ is a Sylow $p$-subgroup of $X$ (\cite{Wolf}, Theorem A).
Thus by a) we get $\Nor_H(P) =P$. \\
c)  This follows immediately from Theorem 1.1 of \cite{GMN} since $H$ is $p$-solvable.
\end{Proof}

We would like to mention here that recently (see \cite{NT}) it has been proved  that for a general finite group  $G$ and $p$ an odd prime the
conditions a) and b) in Proposition \ref{solvable} are equivalent. \\

\medskip

Let $G$ be a $p$-solvable group and let $\varphi\in\IBr_p(G)$. 
We say that a pair  pair $(J,\gamma)$ where $J=J_\varphi \leq G$ and  $\gamma = \gamma_\varphi \in\IBr_p(J)$ is a Huppert pair if 
$\gamma$ is of $p'$-degree and $\gamma^G=\varphi$.
By a result of Huppert (\cite{Navarro}, Theorem 10.11), every $\varphi\in\IBr_p(G)$ has a Huppert pair.

\begin{Proposition} \label{red4} Let $G$ be a $p$-solvable group and let  $\varphi
\in\IBr_p(G)$ with Huppert pair $(J,\gamma)$.
Then the following conditions are equivalent.
\begin{itemize}
\item[\rm a)] $\varphi^\Fix(1)=c_\varphi=\varphi(1)_{p'}.$
\item[\rm b)] $ p\nmid |J:\Kernel\;\gamma|$
and $P^x\cap J\in\Syl_p(J)$ for all $x \in G$.
\end{itemize}
\end{Proposition}
\begin{Proof} Assume that part a) holds. Let $V$ denote the $FG$-module affording $\varphi$ and let $W$ be the $FJ$-module affording $\gamma$.
Since $W^G=V$
 Mackey's Formula and Nakayama's Lemma yield

$$\begin{array}{rcl}
\varphi^\Fix(1) & = & \dim\Hom_{FP}(1_P,V_P)=\dim\Hom_{FG}(1_P^G,V) \\[2ex]
 & = & \dim\Hom_{FG}(1_P^G,W^G)=  \dim\Hom_{FJ}(1_P^G|_J,W) \\[2ex]
 &=&\sum_{x\in P\backslash G/J}\dim\Hom_{FJ}(1_{P^x\cap J}^J,W) \\[2ex]
& \leq &
\sum_{x\in P\backslash G/J}(1_{P^x\cap J}^J,\Phi_\gamma)=
\sum_{x \in P\backslash G/J}\Phi_\gamma(1){|J:P^x\cap J|\over|J|} \\[2ex]
& = &
{1\over|J|}\Phi_\gamma(1)\sum_{x \in P\backslash G/J}|J:P^x\cap J| = \\[2ex]
& = &
{1\over|J|}\gamma(1)|J|_p|G:P| 
=
|G:J|_{p'}\gamma(1)=\varphi(1)_{p'}.
\end{array}$$
Thus the condition
$ \varphi^\Fix(1)=\varphi(1)_{p'}$
 is equivalent to
 \begin{equation}\label{condition}
 \dim\Hom_{FJ}(1_{P^x\cap J}^J,W)=\gamma(1){|J|_p\over|P^x\cap J|}\text{ for every }x.
 \end{equation}
 Observe that we may assume $P\cap J\in\text{\Syl}_p(J)$. Therefore (\ref{condition}) implies
 $$\gamma^\Fix(1)=\dim\Hom_{FJ}(1_{P\cap J}^J,W)=\gamma(1),$$
which means that $P\cap J$ acts trivially on $W$. Thus $\text{Ker}\;\gamma \trianglelefteq J$ contains all the $p$-subgroups of $J$. In turn, this implies that   $P^x\cap J\leq\text{Ker}\;\gamma$ for every $x \in G$. It follows
$$\gamma(1)=\dim\Hom_{FJ}(1_{P^x\cap J}^J,W)=\gamma(1){|J|_p\over|P^x\cap J|},$$
hence $|P^x\cap J|=|J|_p$. Thus $P \cap J$ is a Sylow $p$-subgroup of $J$ for every Sylow $p$-subgroup $P$ of $G$.

Conversely, suppose that part b) holds. If $p\nmid |J:\text{Ker}\;\gamma|$ then for every $x \in G$ we have
$$\gamma(1)= \dim\Hom_{FJ}(1_{P^x\cap J}^J,W),$$
and if moreover $|J|_p=|P^x\cap J|$, we obtain (\ref{condition}) which is equivalent to the statement in a).
\end{Proof}

As a consequence we have the next proposition. Recall that if $G$ is $p$-solvable and $K$ is a $p$-complement, a Fong character associated to $\varphi\in\IBr_pG$ is a $\chi\in\text{Irr}(K)$ such that $\chi^G=\Phi_\phi$. By (\cite{Navarro} Theorem 10.13), any $\varphi\in\IBr_p(G)$ has some associated Fong character.

\begin{Proposition}\label{pcomplement1} Assume that $G$ has a normal $p$-complement $K$. Let $\varphi\in\IBr_p(G)$ and let $\chi\in\Irr(K)$ be an associated Fong character. Then $ \varphi^\Fix(1)=\varphi(1)_{p'}$
 if and only if  the inertial group $\I_P(\chi)$ of $\chi$ in $P$ acts trivially on  $K/\Kernel \chi$.
\end{Proposition}
\begin{Proof} 
First observe that the condition on $\I_P(\chi)$ does not depend on the chosen Sylow $p$-subgroup. 
Moreover, $\I_{P^x}(\chi)$ is a Sylow $p$-subgroup of $J=\I_G(\chi)=(J\cap P)K=\I_{P^x}(\chi)K$ for every $x\in G$.
%Note that for $x \in K$ the group $P^x \cap J = \I_{P^x}(\chi) = \I_P(\chi)^x$ is a Sylow $p$-subgroup of $J$ as well.
Furthermore $\I_P(\chi) \Kernel \chi \trianglelefteq J = \I_P(\chi)K$ if and only if
$$ [\I_P(\chi), K] \leq\I_P(\chi)\Kernel\chi\cap K= \Kernel \chi$$ 
which is equivalent to the fact that $\I_P(\chi)$ acts trivially on $K/\Kernel \chi$.

Assume first that $ \varphi^\Fix(1)=\varphi(1)_{p'}$.
By a result of Dade (see \cite{Navarro},  Theorem 8.13), there exists a $\gamma\in\IBr_p(J)$ such that $\chi=\gamma|_K$. In particular, $\gamma$ has $p'$-degree. Thus $(J,\gamma)$ is a Huppert pair for $\varphi$, and by Proposition \ref{red4}, we obtain $p\nmid|J:{\Kernel}\gamma|$. Hence $\I_P(\chi)\leq {\Kernel}\gamma$. Since ${\Kernel}\gamma\cap K={\Kernel}\chi$ and $J=I_P(\chi)K$ we get $ \I_P(\chi){\Kernel}\chi = {\Kernel}\gamma \trianglelefteq J$.

Conversely, suppose that $T=\I_P(\chi){\Kernel}\chi\trianglelefteq J$. Since $J/T\cong K/{\Kernel}\chi$ we may consider the Brauer character $\gamma$ of $J$ obtained from $\chi$ by inflation. It follows again that $(J,\gamma)$ is a Huppert pair for $\varphi$ and since ${\Kernel}\gamma=T$, we get 
$ p \nmid |J : \Kernel \gamma|$. Recall that for every 
$x\in G$ the intersection $P^x \cap J$ is a Sylow $p$-subgroup of $J$.
Thus condition b) in Proposition \ref{red4} is satisfied  from which $ \varphi^\Fix(1)=\varphi(1)_{p'}$ follows.
\end{Proof}

\begin{Corollary}\label{pcomplement2} Assume that $G$ has a normal $p$-complement $K$. Then $1_P^G$ is completely reducible if and only if 
   $\I_P(\chi)$ acts trivially on $K/\Kernel \chi$  for all $\chi\in \Irr(K)$. 
\end{Corollary}
\begin{Proof} Since $G$ has a normal $p$-complement Green's Indecomposability Theorem (\cite{Huppert-Black}, Chap. VII, 16.2) yields that
all irreducible characters of $K$ are Fong characters. Thus the statement follows immediately by Proposition \ref{pcomplement1}
and Proposition \ref{lower}.
\end{Proof}

\begin{Example} {\rm Let
$F=\operatorname{GF}(3)$ and $p=2$. We put
$$ K = \Bigg\{ \left( \begin{array}{ccc} 1 & * & * \\ 0 & 1 & * \\ 0 & 0 & 1 \end{array} \right) \mid * \in F \Bigg\}$$
and 
$$ P = \Bigg\{ \left( \begin{array}{ccc} a & 0 & 0 \\ 0 & b & 0 \\ 0 & 0 & c \end{array} \right) \mid 0 \not= a,b,c \in F \Bigg\}.$$
Thus $ G=KP$ is the Borel subgroup of $\GL(3,3)$ and  $ K \trianglelefteq G$.  One easily checks that $\OO^{2'}(G) = G$ and $\Nor_G(P) = P$.
From the characters of the extraspecial normal subgroup $K$ we can construct all $ \varphi \in \IBr_2(G)$. There is exactly one
$\varphi$ of degree $6$ and all other Brauer characters $\varphi \not= 1$ are nonlinear and have  $2$-power degree.  
We consider the Brauer character $\varphi$ of degree $6$. The two faithful characters, say $\chi$ and $\chi'$ of $K$ of degree $3$ are Fong characters for $\varphi$ and $ \C_P(K) = 1$. On the other hand $|\I_P(\chi)|=4$. Therefore $1_P^G$ is not completely reducible by
Corollary \ref{pcomplement2}. The example shows that the conditions a) - c) do not characterize groups for which $1_P^G$ is completely reducible.}
\end{Example}

\begin{Corollary} Assume that $G$ has a normal $p$-complement $K$. Then 
$1_P^G$ is completely reducible if for every $g\in P$ we have $\C_{L_g}(g)=1$, where
$$L_g =\cap \, \{T\trianglelefteq K\mid g \in\Nor_P(T), \ \mbox{and $g$ acts trivially on $K/T$} \}.$$
\end{Corollary}
\begin{Proof} 
 Let $g\in P$ and $\chi\in \Irr(K)$ with $g\in\I_P(\chi)$. Let $\phi\in\text{Irr}(L_g)$ such that $\chi|_{L_g}$ contains $\phi$ as a constituent. 
 As $C_{L_g}(g)=1$ and $L_g$ has $p'$-order, one easily sees that there is no nontrivial $L_g$-conjugacy class in $L_g$ which is setwise fixed under the action of $g$. Thus if $1\neq\phi \in \Irr(L_g)$, we have $\phi^g\neq\phi$. Since $g$ permutes the irreducible constituents of $\chi|_{L_g}$ and as the number of those is prime to $p$, we deduce that some is fixed. This forces  $\phi=1_N$.  Thus $L_g\leq\text{Ker}\chi$ and $g$ acts trivially on $K/\Kernel \chi$. 
This shows that $\I_P(\chi)$ acts trivially on $K/\Kernel \chi$, and we may apply Corollary \ref{pcomplement2}.
\end{Proof}

\begin{Lemma} Let $G$ be a group with  normal $p$-complement $K$ and 
let $H=O^{p'}(G)$. If
$$L =\cap \, \{T\trianglelefteq K\mid  P \leq \Nor_G(T), \ \mbox{and $P$ acts trivially on $K/T$} \} $$ then
 $\Nor_H(P)=P$ is equivalent to $\C_L(P)=1$. \end{Lemma}
\begin{Proof}  Let $S=H\cap K$,  hence $H=PS$. We shall prove below that 
$L=S$. From this we see that $\Nor_H(P)=P$ is equivalent to 
$\C_L(P)=\C_S(P)=1$.

Observe that  $L\leq S$ is a consequence of the fact that $S$ is normal in 
$G$ and that $G/S\cong
H/S \times K/S  \cong P \times K/S$. Thus $P$ acts trivially on $K/S$.

  In order to prove $ S \leq L$,  note that $P$ acts trivially on $K/L$. 
Therefore $PL$ is normalized by $K$.  Thus $PL$ is normal in $G$. Consequently 
$H\leq PL$ and therefore
 $S=H\cap K\leq PL\cap K=L$ which completes the proof.
\end{Proof}

\section{The case $c_\varphi = 1$} \label{trivial}

In the introduction we mentioned already that the condition $c_\varphi = 1$ for all $\varphi \in \IBr_p(G)$ implies that $G$ is
solvable and all irreducible $p$-Brauer characters have $p$-power degrees. The latter holds true, for instance,
if $G$ has an abelian $p$-complement 
by (\cite{Navarro}, Theorem 10.13),
since in this case all Fong characters are linear,
or if all $p'$-conjugacy classes have  $p$-power size (\cite{Tiep-Willems}, Corollary 4). 

\begin{Lemma}\label{tec} Let $G=\OO_{p',p,p'}(G)$. Then $\varphi(1)$ is  
a power of $p$ for every $\varphi\in\text{IBr}_p(G)$ if and only if the  
following two conditions are satisfied.
\begin{itemize}
\item[\rm a)] All $p'$-sections of $G$ are abelian.

\item[\rm b)]  $G=\C_G(x)\OO_{p',p}(G)$ for all $x\in \OO_{p'}(G)$.
\end{itemize}
\end{Lemma}
\begin{Proof} First suppose that all irreducible $p$-Brauer character  of $G$
have $p$-power degree. 
Since this condition passes on to normal  
subgroups and quotients we obviously have a).

Let $N=\OO_{p',p}(G)$. If $ \phi \in \IBr_p(N)$ then $\I_G(\phi) = G$ by Clifford, since $ p \nmid |G:N|$. 
 By Brauer's Permutation Lemma for irreducible Brauer characters  
(see \cite{Trento} Theorem 7.5) the group $G$ acts trivially  
on the $N$-conjugacy classes of $p'$-elements of $N$. Now take any  
$g\in G$ and $x\in O_{p'}(G)$. As $x^g\in O_{p'}(G)\leq N$ we get  
$x^g=x^n$ for some $n\in N$. Thus  $gn^{-1}\in C_G(x)$ and part b) follows.

Conversely, assume that a) and b) hold. By part a) and (\cite{Navarro}, Theorem 8.30) we get that all $\phi \in \IBr_p(N)$ have
$p$-power degrees.
Furthermore, part b) together with Brauer's Permutation Lemma imply  $\I_G(\phi) =G$ for all $\phi \in \IBr_p(N)$. 
Thus we only have to check that any such $\phi \in \IBr_p(N)$ can be  
extended to $G$. Since $|G/N|$ and $\phi(1)$ are coprime, by (\cite{Navarro}, Theorem 8.23) 
it suffices to  
prove that $p$ does not divide the order of $\det \phi$. But this is obvious since $\det\phi$ is a linear Brauer character.
\end{Proof}
 
The next result contains in particular the statement of Theorem C. 

\begin{Theorem}\label{full} The following conditions are equivalent. 
\begin{itemize}
\item[\rm a)] $c_\varphi=1$ for any $\varphi\in\IBr_p(G)$,
\item[\rm b)] $\varphi(1)$ is a power of $p$ for every $\varphi \in \IBr_p(G)$, and $G$ is solvable if $p=2$.
\item[\rm c)] $G=O_{p,p',p,p',p}(G)$ and  
$L=O_{p,p',p,p'}(G)/O_p(G)$ satisfies {\rm a)} and {\rm b)} in {\rm Lemma \ref{tec}}.
\end{itemize}
\end{Theorem}
\begin{Proof} We first claim that a) and b) are equivalent. Since  $ 1 \leq \varphi^\Fix(1) \leq c_\varphi$
part a) together with Theorem B 
 imply that $G$ is $p$-solvable.

Hence, by the Feit-Thompson Theorem, $G$ is solvable in case 
$p=2$. Furthermore $ 1 = c_\varphi = \varphi(1)_{p'}$ for all $\varphi \in \IBr_p(G)$. 
Thus we have b). On the other hand, if b) holds then by (\cite{Tiep-Willems}, Corollary 2), the group $G$ is  
solvable. Hence $c_\varphi=\varphi(1)_{p'}=1$ for every $\varphi\in\IBr_p(G)$ and we have a).

Now we prove that b) implies c). Recall that as observed before $G$ is solvable. Thus   
 Theorem 13.10 of \cite{ManzWolf} implies $G=O_{p,p',p,p',p}(G)$.  
Since our hypothesis passes on to  sections, 
 $L$ is characterized by a) and b) of Lemma \ref{tec}.

Finally, c) means, using again Lemma \ref{tec}, that all   
irreducible $p$-Brauer characters of $L$ have degree a power of $p$.  
Obviously, the same happens for $T=O_{p,p',p,p'}(G)$ since normal $p$-subgroups are
in the kernel of irreducible representations in characteristic $p$. 
Thus, by (\cite{Navarro} Theorem 8.30), we get b).
\end{Proof}

\begin{Remark}{\rm Whereas in Sections 2 and 3 we have considered only irreducible Brauer characters in $\IBr_p(\Phi_1)$, we are dealing in this section with the whole character set $\IBr_p(G)$. Using a proof along the lines of the proof of Theorem B, one sees that if $c_\varphi=1$ for any $\varphi\in\IBr_p(\Phi_1)$, then $G$ must be $p$-solvable. But any group with a normal $p$-complement shows that the full characterization in Theorem C doesn't work if one  considers only characters in $\IBr_p(\Phi_1)$. However, it might be true that there were some restriction on the $p$-length.}
\end{Remark}

\begin{Question} Let $G$ be a group with $c_\varphi=1$ for any $\varphi\in\IBr_p(\Phi_1)$. As remarked above $G$ is $p$-solvable. Is there some bound on the $p$-length of such a group $G$?
\end{Question}

\section{The case $c_\varphi = \varphi(1)$}
 Not much is known about $c_\varphi$ for groups which are not $p$-solvable. For example, if $c_1=1$ then $c_\varphi \leq
\varphi(1)$, but this does not hold true in general.
%for all $\varphi \in \IBr_p(G)$.
However, in the average it might be true
as many examples show. To be more precise,
let $\IBr_p(G) = \{\varphi_1, \ldots, \varphi_k\}$. We denote by
$\overline{\varphi}$ the dimension vector
$(\varphi_1(1), \ldots, \varphi_k(1))$ and by $c$ the vector
$(c_{\varphi_1}, \ldots, c_{\varphi_k})$.

Furthermore, let $< \cdot \,, \cdot >$ denote the euclidean scalar
product on $\R^k$ with corresponding norm $\parallel \cdot \parallel$.
In particular, $\parallel \overline{\varphi} \parallel^2 =
\sum_{\varphi \in \IBr_p(G)} \varphi(1)^2.$
With this notation we state

\begin{Conjecture} \label{conj}  The inequality
% $$ \frac{1}{(|G|_p)^2} \sum_{\varphi \in \IBr_p(G)}
% \Phi_\varphi(1)^2 \leq  \sum_{\varphi \in \IBr_p(G)} \varphi(1)^2.$$
$$ \parallel {c} \parallel \ \leq \ \parallel \overline{\varphi}
\parallel $$ holds  true for every finite group.
%with equality if and only if the Sylow $p$-subgroup is normal.
\end{Conjecture}

Based on many examples we conjectured in \cite{Willems} the following.

\begin{Conjecture} \label{conj1} We always have
%\begin{equation} \label{Brauer-degrees} |G|_{p'} \leq \sum_{\varphi
%\in \IBr_p(G)} \varphi(1)^2 \end{equation}$
\begin{equation}   |G|_{p'} \leq \ \parallel \overline{\varphi} \parallel^2 \end{equation}
with equality if and only if the Sylow $p$-subgroup is normal.
\end{Conjecture}

\begin{Proposition} \label{main2} If Conjecture {\rm \ref{conj}} has an
affirmative answer then Conjecture {\rm \ref{conj1}} has an
affirmative answer 
including the characterization of equality.
\end{Proposition}
\begin{Proof}
Conjecture \ref{conj} states that
$$ < c,c> \ \leq \ <\overline{\varphi}, \overline{\varphi} >.$$
Applying the Cauchy-Schwarz inequality we get
\begin{equation} \label{CS}
  <c,\overline{\varphi}>^2 \ \leq \ <c,c> \, <\overline{\varphi},
\overline{\varphi}> \ \leq \  <\overline{\varphi},
\overline{\varphi}>^2.
\end{equation}
This proves  that
$$ (|G|_{p'})^2 = \frac{1}{(|G|_p)^2} (\sum_{i=1}^k \Phi_{\varphi_i}(1)
\varphi_i(1))^2 \ \leq \ (\sum_{i=1}^k \varphi_i(1)^2)^2 = \
\parallel \overline{\varphi} \parallel^4,$$
which is the first part of Conjecture \ref{conj1}.

Next we characterize equality.
Clearly, if the Sylow $p$-subgroup is normal then $p \nmid
\varphi_i(1)$ for all $i$ and by Fong's Theorem,
$ \Phi_{\varphi_i}(1) = |G|_p\,\varphi_i(1)$. Thus we get
$$ |G|_{p'} = \frac{1}{|G|_p} (\sum_{i=1}^k \Phi_{\varphi_i}(1) \varphi_i(1)) =
\sum_{i=1}^k \varphi_i(1)^2 = \  \parallel \overline{\varphi}
\parallel^2.$$

Finally suppose that we have equality in  Conjecture \ref{conj1} which
means that
$$ (|G|_{p'})^2 = \ <c,\overline{\varphi}>^2 \  = \
<\overline{\varphi}, \overline{\varphi}>^2.$$
Thus the inequality in the Cauchy-Schwarz estimation (\ref{CS}) is
an equality. This forces
$$ c = s \, \overline{\varphi} $$
for some $ s \in \N$.
If we define $\psi(g) = |G|_p$ for $p'$-elements $ g \in G$ and
$\psi(g)=0$ otherwise then $\psi$ is a generalized character,
by (\cite{Trento}, Proposition 15.9). Hence $\psi$ is a $\Z$-linear
combination of the characters $\Phi_{\varphi_i}$ according to (\cite{Trento},
Proposition 15.13).
This shows that
  the greatest common divisor of the degrees of $\Phi_{\varphi_i}$ is equal to
$|G|_p$. It follows
  $s=1$, hence $\Phi_{\varphi_i}(1) = |G|_p \, \varphi_i(1)$.
  Thus $$ \Phi_{\varphi_i} = \Phi_1 \, \varphi_i $$
for $i=1, \ldots, k$ and by a result of Brockhaus \cite{Brockhaus},
the Sylow $p$-subgroup of $G$ is normal
which completes the proof.
\end{Proof}

\begin{Remark} {\rm Conjecture \ref{conj1} says that $|G|_{p'} \leq \sum_{\varphi \in \IBr_p(G)} \varphi(1)^2$. The equality 
$$ |G|_{p'} =\sum_{\varphi \in \IBr_p(G)} \varphi(1)^2$$
is equivalent to 
$  c_\varphi = \varphi(1)$ for all $\varphi \in \IBr_p(G) $. This should happen if and only if a Sylow $p$-subgroup is normal in $G$.
 The inequality (\ref{fix}) implies
$$ \sum_{\varphi \in \IBr_p(G)} \varphi(1) \leq \sum_{\varphi \in \IBr_p(G)} \varphi^\Fix(1)\varphi(1) \  \leq \  |G|_{p'}.$$
Thus we have  
$$  \sum_{\varphi \in \IBr_p(G)} \varphi^\Fix(1)\varphi(1) =  |G|_{p'} $$
if and only if $c_\varphi = \varphi^\Fix(1)$ for all $\varphi \in \IBr_p(G) $. According to section \ref{Fix(P)} this happens if and only if
$1_P^G$ is completely reducible. Unfortunately we do not have a full classification in terms of the group structure in this case. Finally,
$$  \sum_{\varphi \in \IBr_p(G)} \varphi(1) =  |G|_{p'} $$
if and only if $c_\varphi =1$ for all $\varphi \in \IBr_p(G) $. A full characterization of groups which satisfy this condition is given in 
section \ref{trivial}.
}
\end{Remark}

\end{document}